\documentclass[a4paper,10pt]{article}
\usepackage{ifstyle}
\usepackage{graphicx}
\usepackage{epsfig} 

\usepackage[english]{babel}               
\usepackage[T1]{fontenc}
\usepackage[latin1]{inputenc}
\usepackage{array}                 
\def\Xint#1{\mathchoice
{\XXint\displaystyle\textstyle{#1}}%
{\XXint\textstyle\scriptstyle{#1}}%
{\XXint\scriptstyle\scriptscriptstyle{#1}}%
{\XXint\scriptscriptstyle\scriptscriptstyle{#1}}%
\!\int}
\def\XXint#1#2#3{{\setbox0=\hbox{$#1{#2#3}{\int}$}
\vcenter{\hbox{$#2#3$}}\kern-.5\wd0}}

\def\dashint{\Xint-}
\newtheorem{theorem}{Theorem}

\newtheorem{lem}{Lemma}

\def\opn#1#2{\def#1{\operatorname{#2} } } 
\opn\Rm{Rm}
\opn\Ric{Ric} 
\opn\Rc{Rc}
\opn\Scal{Sc}
\opn\Tr{Tr}
\opn\Trac{Trace} 
\opn\det{det} 
\opn\div{div}
\opn\Ker{Ker} 
\opn\exp{exp}
\opn\exph{exph}
\opn\Herm{Herm}
\opn\End{End} 
\opn\Hess{Hess}

\newcommand{\R}{\mathbb R}

\newcommand{\E}{\mathbb E}

\newcommand{\contract}{\mathrel{\kern-1.5pt\vrule width6.0pt height0.4pt depth0pt
                \vrule width0.4pt height4.0pt depth0pt}}
\newcommand{\retract}{\mathrel{\kern-1.5pt\vrule width0.4pt height4.0pt depth0pt
                          \vrule width6.0pt height0.4pt depth0pt}}

\begin{document}
\begin{center} 
\Large{\bf{A consequence of a lower bound of the K-energy }}
\\
\vspace{0.4cm}
\large{Nefton Pali}
\end{center} 
\begin{abstract}We prove that over a Fano manifold with K-energy of the canonical class $2\pi c_1$ bounded from below, the Chen-Tian energy functional $E_1$ of the canonical class is allso bounded from below.
\end{abstract}
\section{Introduction}
One of the central philosophical points of modern Differential Geometry is in the idea that the existence of special metrics should imply strong geometric conditions on the manifold. An important example of special metrics are the extremal metrics of a K\"ahler manifold, introduced by Calabi \cite{Cal}. By definition extremal K\"ahler metrics of a positive $(1,1)$-cohomology class are K\"ahler metrics with minimal $L^2$-norm of the scalar curvature with respect to all K\"ahler metrics in the same $(1,1)$-cohomology class. Calabi proved (see \cite{Cal}) that extremal metrics are exactly the K\"ahler metrics with holomorphic gradient of the scalar curvature. To any positive $(1,1)$-cohomology class Calabi and Futaki \cite{Fut} associated an invariant such that its vanishing implies that the extremal metrics of the given class are exactly the metrics with constant scalar curvature. Moreover the vanishing of the Calabi-Futaki invariant is a necessary condition for the existence of metrics with constant scalar curvature.
In the case of Fano manifolds, a positive multiple of the first Chern class admits a metric with constant scalar curvature if and only if the K\"ahler metric is Einstein. The problem of existence of K\"ahler-Einstein metrics over Fano manifolds has been the subject of intense study over the last two decades. Inspired from the work of Donaldson, Mabuchi introduced the K-energy functional \cite{Mab}. In a joint work with Bando, \cite{Ba-Ma} they proved that the existence of a K\"ahler-Einstein metric implies a lower bound of the  K-energy and the uniqueness of the K\"ahler-Einstein metric modulo the action of automorphisms of the manifold. Recently it has been shown by Chen and Tian \cite{Ch-Ti3} that the existence of an extremal metric in a given class implies the lower boundedness of the K-energy and the uniqueness of the extremal metric modulo the action of automorphisms. In the special case of a multiple of the first Chern class we have a very powerfull tool which is the K\"ahler-Ricci flow. Using a similar computation as in \cite{Yau} it can be proved (see \cite{Cao}) that the K\"ahler-Ricci flow admit always a solution for all times. In \cite{Cao} Cao proved that the solution converges to a K\"ahler-Einstein metric if the first Chern class is non positive. In this way he re-proved the celebrated Calabi-Yau theorem \cite{Yau}. In the case of Fano manifolds the K\"ahler-Ricci flow does not allways converge, in fact there are Fano manifolds which do not admit K\"ahler-Einstein metrics \cite{Fut}, \cite{Tia}. In \cite{Ch-Ti1} Chen and Tian give a proof which shows that
the existence of a K\"ahler-Einstein metric over a Fano manifold with nonnegative bisectional curvature implies the convergence of the K\"ahler-Ricci flow. Their hypothesis on the nonnegativity of the bisectional curvature is just a temporary technical assumption. What they really need is that the positivity of the Ricci tensor is preserved under the K\"ahler-Ricci flow. In order to prove the convergence Chen and Tian introduced in \cite{Ch-Ti2} a family of energy functionals $E_k,\,k=0,...,n-1$. Until now, only $E_0$ and $E_1$ plays a key role in the convergence of the K\"ahler-Ricci flow. The functional $E_0$ is just the K-energy. One of the crucial points of their proof is that, over K\"ahler-Einstein manifolds admiting a metric $\omega_0$ such that the nonegativity of the Ricci tensor is preserved under the K\"ahler-Ricci flow with initial metric $\omega_0$, the energy functional $E_1$ is lower bounded from below (so bounded) under the flow. In fact the energy functional $E_1$ is always decreasing along the flows which preserve a certain uniform lower bound of the Ricci tensors of the evolved metrics. Very recently Song and Weinkove \cite{So-We} have answered a question of Chen \cite{Che}, showing that over K\"ahler-Einstein manifolds the energy functional $E_1$ is bounded on the full space of potentials. We have the following general result.
\begin{theorem}
Let $X$ be a Fano manifold such that the K-energy functional of the canonical class $2\pi c_1(X)$ is  bounded from below. Then the Chen-Tian energy functional $E_1$ of the canonical class $2\pi c_1(X)$ is allso bounded from below.
\end{theorem}
This is an immediate consequence of the following theorem.
\begin{theorem}
Let $(X,\omega)$ be a polarised Fano manifold with $\omega\in 2\pi c_1(X)$ and let $\nu_{\omega}$ and $E_{1,\omega}$ be respectively the the K-energy  and the Chen-Tian energy functionals with reference metric $\omega$. Then there exist a constant $C_{\omega}$ depending only on $\omega$ such that for any K\"ahler-Ricci flow $(\omega _t)_{t\in [0,+\infty)}\subset {2\pi c_1}(X)$ we have the inequality
$E_{1,\omega}(\omega_t)\geq 2\nu_{\omega}(\omega_t)+C_{\omega}$ for all $t\in [0,+\infty)$.
\end{theorem}
We consider it necessary to write this note because our proof of theorem 2 is drastically simple. Moreover there are indications that our result will be usefull for a new existence criteria of K\"ahler-Einstein metrics.
\\
\\
{\bf Acknowledgments.} The result presented in this note has been reported during the visit of the author in Princeton University. The author is very grateful to Professors Jean-Pierre Demailly, Joseph Kohn and Gang Tian who made possible his long visit in this institution. The author is especially grateful to Professor Gang Tian for bringing to his attention problems related with the convergence of the K\"ahler-Ricci flow. The author is allso very grateful to the referees for their nice suggestions which have contribute to improve the original version of this note.
	
\section{Energy functionals} 
Let $X$ be a Fano manifold of complex dimension $n$ and let $\omega \in 2\pi c_1(X)$. Let
$
{\cal K}_{2\pi c_1} := \{\omega>0,\,\omega \in 2\pi c_1(X)  \},
$
be the space of K\"ahler metrics in the class $2\pi c_1(X)$.
We consider also the corresponding space of potentials
$
{\cal P}_{\omega}:= \{\varphi \in {\cal E}(X,\R)\,|\,i\partial\bar{\partial}\varphi >-\omega\} 
$
and we define
$\omega _{\varphi}:=\omega +i\partial\bar{\partial}\varphi$ for every $\varphi \in {\cal P}_{\omega}$. 
We will use also the notation $\dashint_X:=(\int_X\omega ^n)^{-1}\int_X$ for the average operator. We remind that the scalar curvature $\Scal(\omega)\in {\cal E}(X,\R)$ of $\omega$ is defined by the formula
$$
\Scal(\omega):=\Trac_{_{\omega} }(\Ric(\omega ))= \frac{2n\Ric(\omega )\wedge \omega ^{n-1} }{\omega ^n},  
$$
where $\Ric(\omega )$ is the Ricci form. This definition coincides with the usual definition of scalar curvature of the Riemannian metric $g=\omega (\cdot, J\cdot)$.
We define the Laplacian  of a function $f$ by the formula
$$
\Delta _{_{\omega}}f:=\mbox{Trace}_{_{\omega} } (i\partial\bar{\partial}f)=
\frac{ 2n\, i\partial_{_J }\bar{\partial}_{_J }f\wedge \omega ^{n-1} }{\omega ^n}.     
$$
Our Laplacian differs by a minus sign from the usual Laplace-Beltrami operator associated to the Riemannian metric $g=\omega (\cdot, J\cdot)$. We remind also that an energy functional $\E$ is a continuous map $\E:{\cal K}_{2\pi c_1}\times{\cal K}_{2\pi c_1}\rightarrow \R$, which satisfies the conditions
\begin{eqnarray*}
\E(\omega_1,\omega_3)&=& \E(\omega_1,\omega_2)+\E(\omega_2,\omega_3)\label{CocicleEner}
\\
\E(\omega_1,\omega_1)&=& 0
\end{eqnarray*}
for every $\omega_1,\,\omega_2,\,\omega_3\in {\cal K}_{2\pi c_1}$.
\\
\\
{\bf The generalized energy functional.}
\\
The generalized energy functional
$J_{\omega}:{\cal P}_{\omega}\rightarrow [0,+\infty)$ is defined by the formula
\begin{eqnarray*}
J_{\omega}(\varphi )&:=&\sum_{k=0}^{n-1} \frac{k+1}{n+1}\,
\dashint\limits_X i\partial \varphi \wedge \bar{\partial}\varphi \wedge\omega ^k\wedge \omega^{n-k-1}_{\varphi}
\\
\\
&=&\dashint\limits_X\varphi \,\omega ^n-\frac{1}{n+1}\sum_{k=0}^n\,\dashint\limits_X\varphi \,\omega ^k\wedge \omega ^{n-k}_{\varphi}.
\end{eqnarray*}
If $(\varphi _t)_{t\in (-\varepsilon ,\varepsilon )} \subset{\cal P}_{\omega}$ is a $\ci$ patht then we have the important formula
\begin{eqnarray}\label{der-JFunc}
\frac{d}{dt}\sum_{k=0}^n\,\dashint\limits_X\varphi_t \,\omega ^k\wedge \omega ^{n-k}_t =(n+1)\dashint\limits_X\dot{\varphi}_t\,\omega ^n_t,
\end{eqnarray}
where $\dot{\varphi}_t:=\frac{\partial}{\partial t}\varphi _t$ and $\omega_t:=\omega_{\varphi_t}$. In fact consider the equalities
\begin{eqnarray*}
&&\frac{d}{dt}\sum_{k=0}^n\,\dashint\limits_X\varphi _t \,\omega ^k\wedge \omega ^{n-k}_t =
\\
\\
&=&
\sum_{k=0}^n\,\dashint\limits_X\left(\dot{\varphi}_t\,\omega ^k\wedge \omega ^{n-k}_t
+
(n-k)\,\varphi _t\,i\partial\bar{\partial}\dot{\varphi}_t\wedge\omega ^k\wedge \omega ^{n-k-1}_t \right)
\\
\\
&=&
\sum_{k=0}^n\,\dashint\limits_X\left(\dot{\varphi}_t\,\omega ^k\wedge \omega ^{n-k}_t
+
(n-k)\,\dot{\varphi}_t\,(\omega_t-\omega )\wedge\omega ^k\wedge \omega ^{n-k-1}_t \right)
\\
\\
&=&
\sum_{k=0}^n\,(n-k+1)\,\dashint\limits_X\dot{\varphi}_t\,\omega ^k\wedge \omega ^{n-k}_t
-
\sum_{k=0}^n\,(n-k)\,\dashint\limits_X\dot{\varphi}_t\,\omega ^{k+1} \wedge \omega ^{n-k-1}_t
\\
\\
&=&
(n+1)\,\dashint\limits_X\dot{\varphi}_t\,\omega ^n_t.
\end{eqnarray*}
If we set $J(\omega,\omega_{\varphi}):=J_{\omega}(\varphi)$ for any $\omega\in {\cal K}_{2\pi c_1}$ then the formula \eqref{der-JFunc} implies that $J$ is an energy functional.
\\
\\
{\bf The K-energy functional of the canonical class $2\pi c_1$.}
\\
Consider the K\"ahler-Ricci flow $(\omega _t)_t,$
\begin{eqnarray}\label{KRic-Flow}  
\frac{d}{dt}\omega _t=\omega _t-\Ric(\omega _t)
\end{eqnarray} 
with initial metric $\omega _0 \in {\cal K}_{2\pi c_1}$. It was proved in \cite{Cao} that the K\"ahler-Ricci flow $(\omega _t)_t$ exists for all $t\in [0,+\infty)$ and $(\omega _t)_t\subset {\cal K}_{2\pi c_1}$. This is because to solve the equation \eqref{KRic-Flow} is sufficient to solve the equation in terms of potentials 
\begin{eqnarray}\label{KRic-FlowPot}
\dot{\varphi}_t=\log\,\frac{\omega _t ^n}{\omega ^n}+\varphi_t-h_{\omega},
\end{eqnarray} 
where $\varphi _t\in {\cal P}_{\omega},\,\omega _t=\omega +i\partial\bar{\partial}\varphi _t$ and $h_{\omega}\in {\cal E}(X,\R)$ is the the real Smooth function defined by the conditions
$$
\Ric(\omega)=\omega +i\partial\bar{\partial}h_{\omega},\quad
\int\limits_X(e^{h_{\omega}} -1)\omega ^n=0.
$$
We remark that to find $\varphi\in {\cal P}_{\omega}$ solution of Einstein the equation 
$
\Ric(\omega _{\varphi} )=\omega _{\varphi}, 
$
is equivalent to solve the equation 
$$
0=\log\,\frac{\omega _{\varphi}^n}{\omega ^n}+\varphi-h_{\omega},
$$
which is also equivalent to the constant scalar curvature equation $\Scal(\omega _{\varphi})=2n$.
This last equation is the  Euler-Lagrange equation of the K-energy functional $\nu_{\omega}:{\cal P}_{\omega}\rightarrow \R$
$$
\nu _{\omega}(\varphi):=
\dashint\limits_X\left(\log\,\frac{\omega _{\varphi} ^n}{\omega ^n}+\varphi-h_{\omega}  \right)\omega^n_{\varphi}
-
\frac{1}{n+1}\sum_{k=0}^n\,\dashint\limits_X\varphi \,\omega ^k\wedge \omega ^{n-k}_{\varphi}
+ \dashint\limits_X h_{\omega}\omega ^n.
$$
In fact for every $\ci$ path $(\varphi _t)_{t\in (-\varepsilon ,\varepsilon )} \subset{\cal P}_{\omega}$ we have the identity
\begin{eqnarray}\label{der-Kenerg} 
\frac{d}{dt}\nu_{\omega}(\varphi_t)=-\frac{1}{2}  \,\dashint\limits_X
\dot{\varphi}_t\Big( \Scal(\omega _t )-2n\Big)\omega ^n_t. 
\end{eqnarray}  
We prove now this identity. Set $\Delta_t:=\Delta_{\omega _t}$ and 
consider the derivative
\begin{eqnarray*}
&&\frac{d}{dt}\,\dashint\limits_X\left(\log\,\frac{\omega _t ^n}{\omega ^n}-h_{\omega} \right)\omega^n_t=
\\
\\
&=&
\dashint\limits_X\left(\frac{d}{dt}\log\,\frac{\omega _t ^n}{\omega ^n}\right)\omega^n_t
+
n\dashint\limits_X\left(\log\,\frac{\omega _t ^n}{\omega ^n} -h_{\omega}\right)\, i\partial\bar{\partial}\dot{\varphi}_t\wedge\omega^{n-1} _t
\\
\\
&=&2^{-1} \dashint\limits_X\Delta_t\dot{\varphi}_t\,\omega^n_t
+
n\dashint\limits_X\dot{\varphi}_t\,i\partial\bar{\partial}\left(\log\,\frac{\omega _t ^n}{\omega ^n}-h_{\omega} \right)\wedge\omega^{n-1} _t
\\
\\
&=&
n\dashint\limits_X\dot{\varphi}_t\,i\partial\bar{\partial}\left(\log\,\frac{\omega _t ^n}{\omega ^n}-h_{\omega} \right)\wedge\omega^{n-1} _t
\\
\\
&=&
-\dashint\limits_X
\dot{\varphi}_t\,2^{-1} \Scal(\omega _t )\,\omega ^n_t
+
n\dashint\limits_X\dot{\varphi}_t\,\omega\wedge\omega^{n-1}_t.
\end{eqnarray*}
Moreover
\begin{eqnarray}\label{DerFun1}
\frac{d}{dt}\,\dashint\limits_X\varphi_t \,\omega^n_t
&=&
\dashint\limits_X\dot{\varphi}_t  \,\omega^n_t+
n\dashint\limits_X\varphi_t \,i\partial\bar{\partial}\dot{\varphi}_t\wedge\omega^{n-1} _t\nonumber
\\\nonumber
\\
&=&\dashint\limits_X\dot{\varphi}_t  \,\omega^n_t
+\dashint\limits_X\dot{\varphi}_t \,(\omega _t-\omega)\wedge\omega^{n-1}_t\nonumber
\\\nonumber
\\
&=&(n+1)\dashint\limits_X\dot{\varphi}_t  \,\omega^n_t
-n\dashint\limits_X\dot{\varphi}_t\,\omega\wedge\omega^{n-1}_t.
\end{eqnarray}
Combining this two equalities with the identity \eqref{der-JFunc} we obtain the identity 
\eqref{der-Kenerg}.
If we set $\nu (\omega ,\omega_{\varphi}):=\nu _{\omega}(\varphi)$ for every $\omega
\in {\cal K}_{2\pi c_1}$ then the formula \eqref{der-Kenerg} implies that $\nu$ is an energy functional. 
We remark that under the  K\"ahler-Ricci flow we have the identity
$$
\Scal(\omega _t)=2n-\Delta _{\omega_t}\dot{\varphi}_t.
$$
Then using the identity \eqref{der-Kenerg} we deduce the inequality
$$
\frac{d}{dt}\nu _{\omega}(\varphi _t)=2^{-1} \dashint\limits_X\dot{\varphi}_t\Delta _{\omega_t}\dot{\varphi}_t\,\omega ^n_t=-n\dashint\limits_X
i\partial \dot{\varphi}_t\wedge \bar{\partial}\dot{\varphi}_t\wedge \omega ^{n-1}_t \leq 0,
$$
which shows that the K-energy decreases under the K\"ahler-Ricci flow. We remind now that the Futaki invariant 
$f_{2\pi c_1}:H^0(X,T_{_{X,J}})\rightarrow \R$ of the K\"ahler class $2\pi c_1(X)$ is defined by the formula
$$
f_{2\pi c_1}(\xi):=\dashint\limits_X\xi .\,h_{\omega }\,\omega ^n
$$
and the definition is independent of the choice of $\omega \in {\cal K}_{2\pi c_1}$. We set by $\mbox{Aut}^0_{_{J}}(X)$ the identity component of the group of $J$-holomorphic automorphisms of $X$. We have the obvious action $\mbox{Aut}^0_{_{J}}(X)\times{\cal K}_{2\pi c_1}\rightarrow {\cal K}_{2\pi c_1}$ given by the pull back.
We remind the following well known fact
\begin{lem}Let $(X,\omega )$ be a compact K\"ahler manifold with $\omega \in {\cal K}_{2\pi c_1}$.
\\
1$)$ If the K-energy $\nu _{\omega}$ is bounded from below then the Futaki invariant $f_{2\pi c_1}$ is zero.
\\
2$)$ If the Futaki invariant $f_{2\pi c_1}$ is zero then the K-energy $\nu _{\omega}$ is $\mbox{Aut}^0_{_{J}}(X)$-invariant.
\end{lem}  
It will be usefull to write the K-energy functional under the following synthetic form
$$
\nu _{\omega}(\varphi):=
\dashint\limits_X\left(\log\,\frac{\omega _{\varphi} ^n}{\omega ^n}-h_{\omega}  \right)\omega^n_{\varphi}
+
\sum_{k=0}^n\frac{a_k}{n+1}\,\dashint\limits_X\varphi \,\omega ^k\wedge \omega ^{n-k}_{\varphi}+C_{0,\omega},
$$
where $a_0=n,\,a_k=-1$ for $k\geq 1$ and $C_{0,\omega}:=\dashint\limits_X h_{\omega}\omega ^n$.
\\
\\
{\bf The Chen-Tian energy functional $E_1$ of the canonical class $2\pi c_1$.}
\\
We define the energy functional $E_{1,\omega}:{\cal P}_{\omega}\rightarrow \R$ of the canonical class $2\pi c_1$ by the formula
\begin{eqnarray*}
E_{1,\omega}(\varphi )&:=&\dashint\limits_X
\Big(\log\,\frac{\omega _{\varphi} ^n}{\omega ^n}-h_{\omega}  \Big)
\Big(\Ric(\omega _{\varphi} )+\omega \Big)\wedge \omega^{n-1}_{\varphi}
\\
\\
&+&
\sum_{k=0}^n\frac{b_k}{n+1}\,\dashint\limits_X\varphi \,\omega ^k\wedge \omega ^{n-k}_{\varphi}+C_{1,\omega},
\end{eqnarray*}
where $b_0=b_1=n-1,\,b_k=-2$ for $k\geq 2$ and 
$$
C_{1,\omega}:=\dashint\limits_X h_{\omega}(\Ric(\omega )+\omega )\wedge \omega ^{n-1}.
$$
This definition coincides with the Chen-Tian definition of the energy functional $E_{1,\omega}$ of the canonical class $2\pi c_1$, (see \cite{Ch-Ti2}). In fact all that we need to show is that for every $\ci$ path $(\varphi _t)_{t\in (-\varepsilon ,\varepsilon )} \subset{\cal P}_{\omega}$ we have the identity
\begin{eqnarray}\label{DerJ1Fun}
\frac{d}{dt}\sum_{k=0}^n\frac{b_k}{n+1}\,\dashint\limits_X\varphi_t \,\omega ^k\wedge \omega ^{n-k}_t=
(n-1)\dashint\limits_X\dot{\varphi}_t(\omega^2_t-\omega^2)\wedge \omega^{n-2}_t.
\end{eqnarray}
Let prove this identity. The formulas \eqref{der-JFunc} and \eqref{DerFun1} implies the equality
\begin{eqnarray}\label{DerFun2}
\frac{d}{dt}\left[\dashint\limits_X\varphi_t\omega^n_t
-
\sum_{k=0}^n\,\dashint\limits_X\varphi_t \,\omega ^k\wedge \omega ^{n-k}_t\right]
=
-n\dashint\limits_X\dot{\varphi}_t\,\omega\wedge \omega^{n-1}_t.
\end{eqnarray}
Moreover we have the equality
\begin{eqnarray*}
\frac{d}{dt}\dashint\limits_X\varphi_t\,\omega\wedge \omega^{n-1}_t
&=&
\dashint\limits_X\dot{\varphi}_t\,\omega\wedge \omega^{n-1}_t
+
(n-1)\dashint\limits_X\varphi_t\, i\partial\bar{\partial}\dot{\varphi}_t \wedge \omega\wedge \omega^{n-2}_t\nonumber
\\\nonumber
\\
&=&
\dashint\limits_X\dot{\varphi}_t\,\omega\wedge \omega^{n-1}_t
+
(n-1)\dashint\limits_X\dot{\varphi}_t \,(\omega_t-\omega)\wedge \omega\wedge \omega^{n-2}_t\nonumber
\\\nonumber
\\
&=&
n\dashint\limits_X\dot{\varphi}_t\,\omega\wedge \omega^{n-1}_t
-
(n-1)\dashint\limits_X\dot{\varphi}_t\,\omega^2\wedge \omega^{n-2}_t.
\end{eqnarray*}
Combining this last equality with the equality \eqref{DerFun2}, we find the identity
$$
-\frac{d}{dt}\sum_{k=2}^n\,\dashint\limits_X\varphi_t \,\omega ^k\wedge \omega ^{n-k}_t
=
-(n-1)\dashint\limits_X\dot{\varphi}_t\,\omega^2\wedge \omega^{n-2}_t.
$$ 
Then using the identity \eqref{der-JFunc} we find that the derivative of the quantity
\begin{eqnarray*}
&\displaystyle{\frac{n-1}{n+1}\sum_{k=0}^n\,\dashint\limits_X\varphi_t \,\omega ^k\wedge \omega ^{n-k}_t
-\sum_{k=2}^n\,\dashint\limits_X\varphi_t \,\omega ^k\wedge \omega ^{n-k}_t
=}&
\end{eqnarray*}
\begin{eqnarray*}
&\displaystyle{=\frac{n-1}{n+1}\left[\dashint\limits_X\varphi_t\,\omega^n_t
+
\dashint\limits_X\varphi_t \,\omega \wedge \omega ^{n-1}_t\right]
-
\frac{2}{n+1}\sum_{k=2}^n\,\dashint\limits_X\varphi_t \,\omega ^k\wedge \omega ^{n-k}_t
=}&
\\
\\
&\displaystyle{=
\sum_{k=0}^n\frac{b_k}{n+1}\,\dashint\limits_X\varphi_t \,\omega ^k\wedge \omega ^{n-k}_t,}&
\end{eqnarray*}
give us the required formula \eqref{DerJ1Fun}. We remind also that \cite{Ch-Ti2} for every $\ci$ path $(\varphi _t)_{t\in (-\varepsilon ,\varepsilon )} \subset{\cal P}_{\omega}$ we have the important identity
\begin{eqnarray}
\frac{d}{dt}E_{1,\omega}(\varphi_t )&=&\dashint\limits_X\Delta_t\dot{\varphi}_t\,\Ric(\omega_t)\wedge \omega^{n-1}_t\nonumber
\\\nonumber
\\
&-&
(n-1)
\dashint\limits_X\dot{\varphi}_t\Big(\Ric(\omega _t )^2-\omega^2_t \Big)\wedge \omega^{n-2}_t,
\end{eqnarray}
which shows in particular that $E_1(\omega,\omega_{\varphi}):=E_{1,\omega}(\varphi)$ is an energy functional. 
\\
We are now in position to prove the theorem 2.
\section{Proof of the theorem 2} 
$Proof$. Under the K\"ahler-Ricci flow
$
\dot{\varphi}_t=\log\,\frac{\omega _t ^n}{\omega ^n}+\varphi_t-h_{\omega},
$
we have the following expression for the K-energy functional
$$
\nu _{\omega}(\varphi_t):=
\dashint\limits_X(\dot{\varphi}_t-\varphi_t)\omega^n_t
+
\sum_{k=0}^n\frac{a_k}{n+1}\,\dashint\limits_X\varphi_t \,\omega ^k\wedge \omega ^{n-k}_t+C_{0,\omega}.
$$
The identity $\Ric(\omega_t)=\omega_t-
i\partial \bar{\partial}\dot{\varphi}_t=
\omega+i\partial \bar{\partial}(\varphi_t-\dot{\varphi}_t)$ implies the following expressions for the energy functional $E_1$ under the K\"ahler-Ricci flow
\begin{eqnarray*}
E_{1,\omega}(\varphi_t )&=&
-\dashint\limits_X(\dot{\varphi}_t-\varphi_t)\wedge 
i\partial \bar{\partial}(\dot{\varphi}_t-\varphi_t)\wedge \omega ^{n-1}_t
+
2\dashint\limits_X(\dot{\varphi}_t-\varphi_t)\,\omega
\wedge \omega^{n-1}_t
\\
\\
&+&
\sum_{k=0}^n\frac{b_k}{n+1}\,\dashint\limits_X\varphi_t \,\omega ^k\wedge \omega ^{n-k}_t+C_{1,\omega}
\\
\\
&=&
\dashint\limits_X i\partial (\dot{\varphi}_t-\varphi_t)\wedge 
\bar{\partial}(\dot{\varphi}_t-\varphi_t)\wedge \omega ^{n-1}_t
+
2\dashint\limits_X(\dot{\varphi}_t-\varphi_t)\,\omega^n_t
\\
\\
&-&
2\dashint\limits_X(\dot{\varphi}_t-\varphi_t)\, i\partial \bar{\partial}\varphi_t\wedge\omega^{n-1}_t
+\sum_{k=0}^n\frac{b_k}{n+1}\,\dashint\limits_X\varphi_t \,\omega ^k\wedge \omega ^{n-k}_t+C_{1,\omega}
\\
\\
&=&
\dashint\limits_X i\partial (\dot{\varphi}_t-\varphi_t)\wedge 
\bar{\partial}(\dot{\varphi}_t-\varphi_t)\wedge \omega ^{n-1}_t
+
2\dashint\limits_X(\dot{\varphi}_t-\varphi_t) \,\omega^n_t
\\
\\
&-&
2\dashint\limits_X i\partial \varphi_t \wedge\bar{\partial}\varphi_t\wedge\omega^{n-1}_t
+
2\dashint\limits_X i\partial\varphi_t\wedge
\bar{\partial}\dot{\varphi}_t\wedge\omega^{n-1}_t
\end{eqnarray*}
\begin{eqnarray*}
+
\sum_{k=0}^n\frac{b_k}{n+1}\,\dashint\limits_X\varphi_t \,\omega ^k\wedge \omega ^{n-k}_t+C_{1,\omega}.
\end{eqnarray*}
Expanding the integral
\begin{eqnarray*}
&\displaystyle{\dashint\limits_X i\partial (\dot{\varphi}_t-\varphi_t)\wedge 
\bar{\partial}(\dot{\varphi}_t-\varphi_t)\wedge \omega ^{n-1}_t=\dashint\limits_X i\partial \varphi_t \wedge\bar{\partial}\varphi_t\wedge\omega^{n-1}_t}&
\\
\\
&\displaystyle{-2\dashint\limits_X i\partial\varphi_t\wedge
\bar{\partial}\dot{\varphi}_t\wedge\omega^{n-1}_t
+
\dashint\limits_X i\partial\dot{\varphi}_t\wedge
\bar{\partial}\dot{\varphi}_t\wedge\omega^{n-1}_t},&
\end{eqnarray*}
we find
\begin{eqnarray*}
E_{1,\omega}(\varphi_t )&=&2\dashint\limits_X(\dot{\varphi}_t-\varphi_t) \,\omega^n_t+
\dashint\limits_X i\partial\dot{\varphi}_t\wedge
\bar{\partial}\dot{\varphi}_t\wedge\omega^{n-1}_t
\\
\\
&-&\dashint\limits_X i\partial \varphi_t \wedge\bar{\partial}\varphi_t\wedge\omega^{n-1}_t
+\sum_{k=0}^n\frac{b_k}{n+1}\,\dashint\limits_X\varphi_t \,\omega ^k\wedge \omega ^{n-k}_t+C_{1,\omega}.
\end{eqnarray*}
Then the trivial identity
\begin{eqnarray*}
&\displaystyle{
\frac{n-1}{n+1}\left[\dashint\limits_X\varphi_t\,\omega ^{n}_t
+
\dashint\limits_X\varphi_t\,\omega \wedge\omega ^{n-1}_t\right]
-
\dashint\limits_X i\partial \varphi_t \wedge\bar{\partial}\varphi_t\wedge\omega^{n-1}_t=}&
\\
\\
&\displaystyle{=\frac{1}{n+1}\left[2n\dashint\limits_X\varphi_t\,\omega ^{n}_t
-
2\dashint\limits_X\varphi_t \,\omega \wedge\omega ^{n-1}_t\right]},&
\end{eqnarray*}
implies the remarkable inequality
$$
E_{1,\omega}(\varphi_t )=2\nu _{\omega}(\varphi_t)+\dashint\limits_X i\partial\dot{\varphi}_t\wedge
\bar{\partial}\dot{\varphi}_t\wedge\omega^{n-1}_t+C_{\omega}\geq
2\nu _{\omega}(\varphi_t)+C_{\omega},
$$
where $C_{\omega}$ is a constant depending only on $\omega$.\hfill$\Box$

\vspace{1cm}
\noindent
Nefton Pali
\\
E-mail address: \textit{npali@Math.Princeton.EDU}
\end{document}